\documentclass[reqno]{amsart}

\usepackage[english]{babel}
\usepackage[all]{xy}
\usepackage {amssymb}

\usepackage{times}
\usepackage[T1]{fontenc}
\usepackage[latin1]{inputenc}
\usepackage{graphicx}

\theoremstyle{plain}
  \newtheorem{theorem}{Theorem}[section]
  \newtheorem{corollary}[theorem]{Corollary}
  \newtheorem*{corollary*}{Corollary}
  
  \newtheorem*{lemma*}{Lemma}
  \newtheorem{proposition}[theorem]{Proposition}
\newtheorem{conjecture}[theorem]{Conjecture}

\theoremstyle{definition}
  \newtheorem{definition}[theorem]{Definition}

\DeclareMathOperator{\id}{id}

\DeclareMathOperator{\grad}{grad}

\DeclareMathOperator{\pr}{pr}
\DeclareMathOperator{\ind}{ind}

\begin{document}

\title{
Dynamics, Spectral Geometry and Topology }

\author{Dan Burghelea}

\address{Department of Mathematics,
         The Ohio State University, 
         231 West 18th Avenue,
         Columbus, 
         OH 43210, USA.}

\email{burghele@mps.ohio-state.edu}




\thanks{
        Partially supported by NSF grant  no
       MCS 0915996}


\begin{abstract}
The paper is an informal report 
on joint work with Stefan Haller on  Dynamics in relation with Topology and Spectral Geometry. By  dynamics  one means  a
smooth vector field on a closed smooth manifold;  the elements of dynamics of concern are the rest points, instantons and closed trajectories.  
One discusses their  counting in the case of a generic vector field which has  some additional properties  satisfied by a still very large class of  
vector fields. 
 \end{abstract}

\maketitle {}
 
\setcounter{tocdepth}{1}
\tableofcontents

\section{Introduction}

This paper provides an informal presentation \footnote{ lecture presented at the 	
"Alexandru Myller" Mathematical Seminar 
Centennial Conference, June 21-26, 2010, Iasi, Romania} of results obtained in collaboration with Stefan Haller;  some of them have been already published in \cite {BH03} and 
\cite {BH07}. 
 We follow here an unpublished version \cite {BH04} of the paper \cite {BH03} which provides a slightly different formulation of the  results. 
There is an additional  gain in generality of the results presented in this paper. The rest points of the  vector fields under consideration rather than of "Morse type" are "hyperbolic".

By dynamics we mean  a smooth vector field $X$ on a closed smooth manifold $M.$  The  elements of dynamics are the rest points, instantons and closed trajectories.  In the case of a generic vector field the rest points are hyperbolic  and because $M$ is closed they are finitely many, but the instantons and closed trajectories are at most countable and  possibly infinitely many.
However, under additional hypotheses,  in each homotopy class of continuous paths between two rest points and in each  homotopy class of closed curves there are finitely many instantons respectively finitely many  closed trajectories.  Therefore  the instantons and closed trajectories can  be counted by "counting functions" defined on the set of corresponding homotopy classes.  When $M$ is connected, in the first case this set is in bijective correspondence with the set of elements of the fundamental group and in the second case with the set of conjugacy classes of elements of the fundamental group.  
For a still large class of vector fields the counting functions  referred above  have  Laplace transforms with respect to some cohomology classes  
$\xi \in H^1(M\;\mathbb R)$ and these  Laplace transforms   are holomorphic functions in some regions of the complex plane. These holomorphic functions  can be described using differential  geometry and topology.

In this paper we express these holomorphic functions in terms of geometric data independent of closed trajectories and instantons (Theorems \ref{th3} and \ref{th4} below).
We also  provide a precise description of the class of vector fields for which our results hold and  show that this class  is large enough to insure that the results are relevant (Theorem \ref {th1} below). These results are not true for all vector fields  but provide some patterns which might be recognized in a much larger class of vector fields than the one we describe and they deserve  additional investigations.

Section 2  of this paper provides the definition of the above mentioned elements of dynamics  and of the properties the vector fields are supposed to have, in order to   make  our theory work,  and this section is referred to as Basics in Dynamics. Section 3 discusses the spectral geometry needed to formulate the results stated in Theorems \ref{th3} and \ref{th4} which provide directly or indirectly interpretation of the Laplace transform  of the counting functions of instantons and  closed trajectories. Section 4  treats some relations with topology. This section is very sketchy because more substantial results will need additional lengthy definitions. We review only Novikov results about the rest points and Huchings-Lee and Pajitnov results about closed trajectories. To find more of what we know at this time the  interested reader is invited to consult \cite {BH07}.

\section {Basics in Dynamics}
Let $M$ be a smooth manifold, $TM\overset{p}{\to} M$ the tangent bundle and $X$ a smooth vector field. Recall that a smooth vector field can be regarded as a smooth map $X: M \to TM$ s.t.  $ p \cdot X=id.$ 

Any vector field $X$ has a flow $\Psi_t: M\to M ,$ a smooth one parameter  group of diffeomorphisms. We denote by $\theta_m: \mathbb R\to M,$ the unique trajectory which passes through $m\in M$ at $t=0,$ which is exactly $\theta_m(t):= \Psi_t (m).$ 
\vskip .2in
{\bf Rest points:} The set of rest points is the set $\mathcal X= \{x\in M | X(x)=0\}$ where the vector field $X$ vanishes. 

At any such point the differential of the the map $X:M\to TM$ defines the linear map $$D_x(X): T_xM\to T_x M$$
 with $T_x(M)= p^{-1}(x),$ 
 as follows. Choose an open neighborhood $U$ of $x$ in $M$, and a trivialization of the tangent bundle 
above $U$, $\theta:TU\to U\times T_x(M)$, with $\theta|_{T_x(M)}=\id$. Consider 
$X^\theta:=\pr_2 \circ  \theta \circ X:U\to T_x(M)$ with $\pr_2$ the projection on the second component. 
Clearly $X^\theta (x)=0$. Observe that $D_x(X^\theta)$ is independent of $\theta$, which justifies  the notation $D_x(X):=D_x(X^\theta).$ 

\begin{definition}  A rest point  $x\in \mathcal X$ is called {\bf nondegenerate = hyperbolic} if  the real part  $\Re \lambda\ne 0$ of each eigenvalue 
$\lambda$ of $D_x(X)$ is different from zero. 
For $x$ hyperbolic rest point the number of the eigenvalues (counted with multiplicity) whose real part is larger than $0$ is called the {\bf Morse index} of $x.$ 
\end{definition} 
\vskip .2in
{\bf Stable/Unstable sets:} For a rest point $x$ the set $W_x^\pm := \{y\in M | \lim_{t\to \pm \infty} \theta_x(t)=x\}$ is called stable /unstable set. If $x$ is nondegenerate and has Morse index $q$ then, by Perron-Hadamard theorem \cite {K95},  $W^+_x$ resp. $W^-_x$ is the image of a one to one immersion  
$$i^+ _x: \mathbb R^{(n-q)}\to M,$$
resp. $$i^-_x: \mathbb R^q\to M,$$
and is called the {\bf stable manifold} resp. {\bf unstable  manifold}  of $x.$
\vskip .2in
{\bf Trajectories:} A smooth map  $\theta:\mathbb R\to M$ is a {\bf parametrized trajectory} if 
$$d\theta(t)/ dt= X(\theta(t)).$$
Two trajectories $\theta_1$ and $\theta_2$ are regarded as equivalent ($\theta_1\equiv \theta_2$) if $\theta_1(t)= \theta_2 (t+a)$ for some real number $a.$

A {\bf nonparametrized trajectory} is an equivalence class $[\theta]$ of parametrized trajectories.

\noindent If $m'= \Psi_t (m)$ the differential $D_m(\Psi_t): T_m(M)\to T_{m'}(M)$ induces 
$\hat {D}_m(\Psi_t): T_m(M)/ T_m(\Gamma) \to T_{m'}(M)/ T_{m'}(\Gamma)$
where  $\Gamma$ denotes the  trajectory containing $m$ and $m'$ and is referred to as {\bf Poincar\'e map} along $\Gamma.$ 

An {\bf instanton} between the rest points $x$ and $y$  is an isolated \footnote {Isolated here means that there exists an open neighborhood $U$ of the underlying set of $[\theta]$ which does not contain any other trajectory between these rest points.} 
nonparametrized trajectory $[\theta]$ with 
$\lim _{t\to -\infty}\theta(t)= x$ and $\lim _{t\to \infty}\theta(t)= y.$
The set of instantons from $x$ to $y$ is denoted by $\mathcal I_{x,y}.$

A {\bf closed trajectory} is a pair $\hat {[\theta]} = ([\theta], T)$ s.t. $\theta(t+T)= \theta(t).$  The number $T$ is called the {\bf time period} of  $\hat{[\theta]}$
as opposed to the {\bf period} $p([\hat \theta])$ of $\hat {[\theta]}$ introduced below. 
The set of closed trajectories is denoted by $\mathcal C.$

\vskip .2in
{\bf Nondegeneracy:}
\begin {definition}:

\begin{enumerate}
\item  A rest point $ x\in \mathcal X$ is called nondegenerate if is hyperbolic, in which case it  has a Morse index $\ind(x)\in \mathbb Z_{\geq 0}$.
The set of rest points of index $k$ is denoted by $\mathcal X_k.$
\item An instanton $[\theta]$ from $x$ to $y,$ with $x, y$ nondegenerate rest points is nondegenerate if the maps $i^-_x$ and $i^+_y$ are transversal at any point   of $[\theta] \subset M.$ Equivalently $W^-_x\pitchfork W^+_y$ along $[\theta],$   in which case $\ind(x)- \ind(y)= 1.$  Moreover,  orientations $o_x$ of $W^-_x$ and $o_y$ of $W^-_y$  induce an orientation of $[\theta]$ and implicitly a sign $\epsilon ^{o_x, o_y}([\theta])$
with  $$\epsilon ^{o_x, o_y}([\theta])= +1\  \rm {or}\  -1$$
 if  the induced orientation is consistent or not  with the orientation from $x$ to $y.$

\item A closet trajectory $[\hat \theta]= ([\theta], T)$ is nondegenerate if the Poincar\'e map 
$$\hat {D}_m(\Psi_T): T_m(M)/ T_m(\Gamma)\to T_m(M)/ T_m(\Gamma)$$
induced from the linear map  $ {D}_m(\Psi_T): T_m(M)\to T_m(M)$
for some $m$ (and then for any other $m\in [\theta]$)  satisfies $\det (\hat {D}_m(\Psi_T)- \lambda Id)\ne 0,\  |\lambda|=1$. In this case denote by  $$\epsilon ( [\hat\theta])= \rm {sign} \det (\hat {D}_m(\Psi_T)-Id).$$
 As the closed trajectory $[\hat \theta]$ defines a map  $\underline{\hat{[\theta]}}: S^1\to M$ 
one denotes by $K([\hat{\theta}])$ the set of integers so that $\underline{\hat{[\theta]}}$ factors by a self map of $S^1$ of degree $k.$
Define  the {\bf period} of $[\hat{\theta}]$  by 
$$p (\hat{[\theta]}): = \text{sup}\  K(\hat{\theta}).$$  
 \end{enumerate}
\end{definition}

\begin{definition}: 

A closed one form $\omega \in \Omega^1 (M), $ ( $d(\omega)=0$) is called {\bf Lyapunov} for $X$ if:
\begin{enumerate}
\item  $ \omega(X) \leq 0$ and  
\item $\omega(X)(x)=0$ iff  $X(x)=0.$
\end{enumerate}

A Lyapunov  form $\omega$ which is Morse (i.e.  locally  is the differential of smooth function with non degenerate critical points) is called {\bf Morse--Lyapunov}.
 
The cohomology class $\xi \in H^1(M;\mathbb R)$ is called Lyapunov for $X$ if it contains  Lyapunov forms for $X.$
\end{definition} 

It is straightforward to check that given a Lyapunov form $\omega$ for $X,$ a neighborhood $U$ of  $\mathcal X$  and $r\geq 2,$ one can find a   smooth function $f$  arbitrary small in $C^r$ topology which  vanishes on $\mathcal X,$ has support in $U,$   and  with $\omega' = \omega +d f$   a Morse--Lyapunov form for $X.$  So if a vector field which satisfies {\bf H} (see below  the definition of {\bf H}) has a Lyapunov form  $\omega$ it also has  a Morse--Lyapunov form $\omega'$ arbitrary close to $\omega.$

\vskip .2in
{\bf Properties of a smooth vector field: }
\begin{itemize}
\item {\bf G} (Genericity)= {\bf H}(Hyperbolicity) + {\bf MS}(Morse Smale) +{\bf NCT}(Nondegenerate closed trajectories)

 Property {\bf H} requires all rest points are nondegenerate;
 
 Property {\bf MS} requires that $i^+_x$ and $i^+_y$
 are transversal for any pair of rest points;
 
 Property {\bf NCT} requires that all closed trajectories are nondegenerate.
\item { \bf L}( Lyapunov) 

Property {\bf L} is satisfied if the set of   Lyapunov cohomology classes in $H^1(M;\mathbb R)$ is nonempty, equivalently if $X$ admits a Lyapunov form. 
\end{itemize}

Let $X$ be a vector field  which satisfies property {\bf H} and $g$ a Riemannian metric on $M.$ For any $x\in \mathcal X$ denote by $\rm {Vol} (B_x(r))$ the volume of the ball centered at $0$ of radius $r$ in $\mathbb R^{\ind x}$ w.r. to the pull back  of the Riemannian metric $g$ by $i^-_x.$ 
\begin{itemize}
\item {\bf EG}(Exponential growth)

The  vector field $X$ has  exponential growth property  {\bf EG}  at the  rest point $x$ if for some ( and therefore every)  Riemannian metric $g$ on $M$ there exists $C> 0 $ so that $\rm {Vol (B_x(r))} \leq e^{Cr}$ for all $r\geq 0.$  It has  exponential growth property {\bf EG} if it has  property {\bf EG}  at all rest points. 
\end{itemize}

These properties are satisfied by a large class of vector fields, as the following theorem indicates.  
\begin{theorem}\label{th1}

1. (Kupka--Smale) For any $r$ the set of vector fields which satisfy $\bold G$ is residual in the $C^r$ topology. 

2. (Smale) Suppose $X$ is a smooth vector field which satisfies  {\bf L}. Then in any $C^r$ neighborhood of $X$ there exists vector fields which coincide with $X$ in a neighborhood of $X$ and satisfy $\bold L$ and $\bold G = \bold H+ \bold{MS}+ \bold{NCT}$.  

3. Suppose $X$ is a smooth vector field which satisfies {\bf L}. Then in any $C^0$ neighborhood of $X$ there exists vector fields which coincide with  $X$ in a neighborhood of $X$ and satisfy $\bold L, \bold G, \rm{and}\  \bold {EG}.$

\end{theorem} 
One expects the following Conjecture to be true.
\begin{conjecture} Statement (3) in the theorem above remains true for an arbitrary $r$.
\end{conjecture}
 This is indeed the case for $n=2.$  In fact for the rest points of Morse index $0, \  1$ and $n$  {\bf EG}  holds;  in the case of index $0$ and $n$ one has nothing to verify. 

In \cite{BH03} we have also considered a stronger version of $\bold{EG}$ referred there as $\bold{SEG}.$ Theorem \ref {th1}  remains true for $\bold{SEG}$ replacing $\bold{EG}.$ Property {\bf SEG} is of the same nature but takes a little longer to describe.

 A very readable proof  of Theorem \ref{th1} (1)  is contained in \cite {P67}.
An inspection of the proof of (1) leads easily to (2) which can also be derived from a slightly stronger version of (1).
The proof of (3) is considerably more elaborated. A complete proof is contained in \cite {BH03} and uses the work of Pajitnov \cite {P99}.

In view of the compacity of $M$ 
Property {\bf L} insures that there are only finitely many rest points, Property {\bf MS} insures that there are at most countable number of instantons and Property {\bf NCT} insures that there are at most  countable number of  closed trajectories.\footnote {In fact much more is true: For any positive real number $T$ the set of closed trajectories  with time period smaller that $T$ is finite
 and similarly the set of instantons whose needed time to go from $i^-_x (S^-_x)$ to $i^+_y(S^+_y)$ is smaller than $T$ is finite. Here $S^-_x$ resp. $S^+_y$ denote the unit sphere in $R^{\ind x}$ resp $R^{n-\ind y}.$}  

The following proposition is of crucial importance.

\begin{proposition}\label {P2.7}
Suppose $X$ satisfies $\bold G$ and $\omega$ is Lyapunov for $X$ representing the cohomology class $\xi.$  
Then for 
any real number $R$ one has:

1. 
The set of instantons $[ \theta] $ so that 
$\int _{[\theta]}\omega < R$ is finite. 

2 . The set of closed trajectories $[\hat \theta] $ so that 
$\xi ([\hat \theta] )< R$ is finite. 
\end{proposition}
Statement (1) is due to Novikov  \cite {N93}.
Statement (2) is due to Fried and Hutchings-Lee  \cite {HL99}. 
 
Proposition  \ref {P2.7} indicates that despite their infiniteness,  the instantons and the closed trajectories  can be counted  with the help of counting functions.
To explain this we need additional definitions. 
 
 For $M$ a closed connected smooth manifold and $x,y\in M$ define:  
 \begin{itemize}
 \item
  $\mathcal P_{x,y}: = $homotopy classes of continuous paths from $x,y$. This set can be put in bijective correspondence to the fundamental group of $M.$
 \item
$ [S^1, M]$ the set of homotopy classes of continuous maps from $S^1$ to $M.$ This set is in bijective correspondence with the conjugacy classes of elements of the fundamental group of $M.$ 
\end{itemize} 

Let  $\xi\in H^1(M;\mathbb R)$ and  $\omega$  a  closed one form representing $\xi.$
Define 
\begin{enumerate}
\item $\omega(\alpha)= \int _{\alpha} \omega,\ $ for $\alpha \in \mathcal P_{x,y},$
\item $\xi (\gamma): = \int_\gamma\omega ,\ \ $ for  $\gamma\in [S^1, M].$
\end{enumerate} 

Let $X$ be a smooth vector field on $M$ which satisfies Property {\bf G}. As noticed 
the set $\mathcal X_k$ of rest points of Morse index $k$ is finite. Denote their number by $n_k.$

Suppose in addition that $X$ has Property {\bf L}. 

In view of Proposition \ref{P2.7} for $x\in\mathcal X_{k+1}, y\in \mathcal X_k,$  $o_x$ orientation of $W^-_x$ and $o_y$ orientation of $W^-_y$ define  the counting functions
$\mathcal I^{o_x, o_y}_{x,y}: \mathcal P_{x,y}\to \mathbb Z$ and  
$\mathcal Z: [S^1, M]\to \mathbb Q$ 
by
$$\mathcal I^{o_x, o_y}_{x,y}(\alpha) :=  \sum _{[\theta]\in \alpha} \epsilon^{o_x, o_y} ([\theta])$$
with $[\theta]$ instanton in the homotopy class $\alpha\in \mathcal P_{x,y}$
and 
$$\mathcal Z (\gamma):= \sum _{[\hat \theta]\in \gamma} \epsilon ([\hat\theta])/ p ([\hat\theta])$$
with  $[\hat \theta]$ closed trajectory  in the homotopy class $\gamma\in [S^1, M].$
\vskip .1in
{\bf Laplace transform:}

If $\omega $ is a closed one form representing a cohomology class $\xi$
one  consider  the following formal expressions:
$$ I^ {o_x, o_y,\omega}_{x,y}(z)=: \sum _{\alpha \in \mathcal P_{x,y}}\mathcal I^{o_x, o_y}_{x,y} (\alpha) e^{-z\omega (\alpha)} $$ and 
$$ Z^\xi(z)=: \sum _{\gamma \in [S^1, M]} \mathcal Z (\gamma)e^{-z\xi (\gamma)} $$
and one can ask when they define holomorphic functions in some parts  of the complex plane.

Note that 
if $\omega'= \omega +dh$ then  $ I^{\mathcal O, \omega} _{x,y}(z)= e^{z(h(y)-h(x)} I^{\mathcal O, \omega'} _{x,y}(z),$ and a change of the  orientations $\mathcal O$ might change  the sign of the function  $ I^{\mathcal O,\omega} _{x,y}(z).$ So an affirmative or negative answer to the question above for the second formal expression depends only on the cohomology class of $\omega.$


\begin{theorem} \label {th2} 
Suppose $X$ satisfies  { \bf G} and {\bf L} and  $\mathcal O= \{ o_x, x \in \mathcal X \}$ is a collection of orientations of the unstable manifolds $W^-_x.$

1. If  $I^{\mathcal O,\omega}_{x,y}(z)$  is absolutely convergent for $\Re z>\rho$ then $Z^\xi (z)$ is absolutely convergent for 
$Re z>\rho.$

2. If $X$ satisfies {\bf G},  {\bf L} and {\bf EG} then there exists $\rho \in \mathbb R$ so that $I^{\mathcal O,\omega}_{x,y}(z)$ and $Z^\xi(z)$ are absolutely convergent
for $\Re z>\rho.$
Moreover, for any $u\in \mathcal X_{k+1}$ and  $w\in \mathcal X_{k-1}$ 
$$\sum_{v\in \mathcal X_k} I^{\mathcal O,\omega}_{u,v}(z) \cdot I^{\mathcal O,\omega}_{v,w}(z)=0.$$
\end{theorem}

Note that the formal sums $I^{\mathcal O,\omega}_{x,y}(z)$ and $Z^\xi (z)$ are the Dirichlet series as defined in \cite{S}. The first is associated with the discrete sequence of real numbers $\omega(\alpha), \alpha \in \mathcal P_{x,y}$ and the corresponding numbers $\mathcal I^{o_x, o_y}_{x,y} (\alpha) \in \mathbb Z\subset \mathbb C.$ The second is associated with the discrete sequence $\xi(\gamma), \gamma \in [S^1, M]$ and the corresponding numbers $\mathcal Z(\gamma)\in\mathbb Q\subset \mathbb C.$ Proposition \ref {P2.7} insures that the formal sums $I^{\mathcal O,\omega}_{x,y}(z)$  and $Z^\xi (z)$
 are Dirichlet series and consequently have an abscissa of convergence $\rho \leq \infty .$ They define holomorphic functions provided the abscissa of convergence is $\ne \infty.$ 

Theorem \ref{th2} is the general result of our work . Ultimately to prove it boils down to show that the abscissa of convergence for the Dirichlet series 
$I^{\mathcal O,\omega}_{x,y}(z)$ and $Z^\xi (z)$ are finite. For this purpose it suffices to consider these series for $z$ a real parameter. This will bring us to Witten deformation described in the next section.  

\vskip .2in 
Suppose that $X$ satisfies {\bf G, L} and {\bf EG} and $\mathcal O=\{o_x, x\in \mathcal X\}$ is a collection of orientations for the unstable manifolds $W^-_x, x\in \mathcal X.$

In view of Theorem  \ref {th2} (2) one  can use the functions $ I^{\mathcal O, \omega}_{x,y} (z)\ $ to define  for any $z, \ \Re z > \rho$ a holomorphic family  of cochain complexes  with base 
$\mathbb C^\ast (M, X, \omega, \mathcal O)(z): = (C^\ast (M,X), \delta^\ast _{\mathcal O, \omega}(z))$ with 

\begin {itemize}
\item $ C^k((M,X):= \rm {Maps} (\mathcal X_k, \mathcal C).$
The base is determined by characteristic functions associated to the rest points. 
\item 
$\delta^\ast _{\mathcal O, \omega}(z); C^\ast (M; X)\to C^{\ast +1}(M, X),$ given by $$\delta^k _{\mathcal O, \omega}(z) (f) (u):= \sum _{v\in \mathcal X_{k}}
I^{\omega}_{u,v}(z) f(v), \ \ u\in \mathcal X_{k+1}$$
\end{itemize}
If $\omega_1, \omega_2$  are two Lyapunov forms representing the same cohomology class $\xi$   and $\mathcal O_1, \mathcal O_2$  are 
two sets of orientations then  there exists a canonical isomorphism between the cochain complexes
$(C^\ast (M, X),\delta ^\ast _{\mathcal O_1, \omega_1} (z))$ and $(C^\ast (M, X),\delta ^\ast _{\mathcal O_2, \omega_2} (z)).$

The isomorphism send the base element corresponding to the rest point $u$ into $\pm e^{z h(u)}$  with $\omega_2= \omega_1 + df$ 
with $\pm$ if $o_{1,u}$ is the same or not with $o_{2,u}.$ We can therefore denote this holomorphic family of cochain complexes, well defined up to an holomorphic isomorphism, 
by $\mathbb C^\ast (M, X,\xi)(z).$

\vskip .2in

\section {Spectral Geometry}
In this section we will describe a few results in geometric analysis and use them to express the holomorphic functions $I^{\mathcal O, \omega}_{x,y}(z)$
and $Z^\xi (z)$ in terms of more geometric invariants.

\vskip .1in
{\bf Witten deformation:}  

Let $\omega$ be a real valued closed one form, $\Omega^\ast (M)$ the real or complex valued differential forms  and $d^\ast_\omega(t): \Omega^\ast (M)\to \Omega^{\ast+1} (M)$ the perturbed exterior differential  defined by $$d^\ast _\omega (t):= d + t \omega\wedge$$ with $t\in \mathbb R.$ 
The family of cochain complexes of deRham type, $(\Omega^\ast (M), d^\ast_\omega(t))$  is referred to as the {"Witten deformation"} (of the deRham complex $(\Omega^\ast(M), d^\ast)$).
 
 Suppose that $M$ is endowed with a Riemannian metric $g.$ Then a differential operator from $\Omega^\ast(M)$ to $\Omega^\ast (M),$   in particular for  $d^\ast_\omega (t):\Omega^\ast (M)\to \Omega^{\ast +1}(M),$ has a  formal adjoint
 $(d_\omega(t)^\sharp)^\ast :\Omega^\ast  (M)\to \Omega^{\ast-1}(M).$
 Following Witten \cite {W82}
 consider $\Delta^\omega_\ast  (t)=(d_\omega(t)^\sharp)^{\ast +1} \cdot d_\omega ^\ast  (t) + d_\omega(t) ^{\ast-1}\cdot (d_\omega (t)^\sharp)^\ast$ which is equal to 
 $$ \Delta_\ast  + t(L_X + L_X^\sharp) + t^2 ||X||^2$$  with $X=-\grad _g(\omega),$ $L_X$ the Lie derivative w.r. to $X,$  $L^\sharp_X$  the formal adjoint 
 of $L_X$ and $||X||^2$ the operator of multiplication by the square of the fiber-wise norm of $X.$  

The operators  $\Delta_\ast^\omega(t)$ are  a zero order perturbation of   $\Delta_q,$ the Laplace--Beltrami operators associated to the Riemannian metric $g.$ 

 Recall that a closed one form is called Morse  if locally is the differential of a function with all critical points  are nondegenerate.  
 As already noticed if $X$ satisfies {\bf H} and {\bf L} then it admits Morse Lyapunov form. 

We will consider the Witten deformation for $\omega$ a  closed Morse one form  and a Riemannian metric flat near the rest points
\footnote {We believe that the flatness requirement is not necessary but  considerably more effort is needed to finalize the arguments without this hypothesis.} 
 of $-\grad_g\omega.$
 Let $n_k$ be the number of rest points of Morse index $k.$

\begin {proposition} There exists positive constants $C_1, C_2, C_3, T$ so that for $t\geq T$ exactly $n_k$ eigenvalues of $\Delta^{\omega, \text{sm}}_k(t)$  counted with their multiplicity are smaller that $C_1 e^{-C_2 t}$ and  all others are larger than $C_3 t.$ 
\end{proposition}

This is a known observation  first made by Witten  \cite {W82} for a Morse exact one form but provable by the same arguments for a Morse  closed one form.

As a consequence, for $t$ large enough  say $t> \rho$,   there is a canonical orthogonal decomposition 
$$(\Omega^\ast (M); d^\omega (t)):= (\Omega^\ast(M)_{\text {sm}} (M)(t), d^\omega (t)) \oplus  (\Omega^\ast(M)_{\text {la}} (M)(t), d^\omega (t))$$
which diagonalizes $\Delta^\omega_\ast (t)$
i.e., $\Delta^\omega_\ast (t) = \Delta^{\omega, \text{sm}} _\ast (t)\oplus \Delta^{\omega,\text{la}}_\ast (t)$ for $t>\rho.$  
The small  resp. large  complex is generated by the eigenforms corresponding to the small resp. large  eigenvalues  i.e. the eigenvalues which for $t >\rho$  are  bounded from above resp. below by a fix number, say $1.$  The small complex $(\Omega^\ast(M)_{\text {sm}} (M)(t), d^\omega (t))$ is finite dimensional with $
\dim \Omega^k(M)_{\text {sm}} (M)(t) = n_k$  while the large complex is acyclic. The number $\rho$ can be any number which insures 
that  $t>\rho $ implies    $C_1 e^{-C_2 t} <1< C_3 t.$

\begin{theorem} \label {th3}
Suppose $X$ satisfies  {\bf G} and  {\bf EG}, $\omega$ is Morse Lyapunov  for $X$ and $\mathcal O$  is a collection of orientations of the 
unstable manifolds $W^-_x, x\in \mathcal X.$   Choose a Riemannian metric on $M$ which is flat near the rest points of $X.$ For any $x\in \mathcal X$ let $h_x: \mathbb R^{\ind x}\to \mathbb R$ be the only smooth function which satisfies  $dh_x= (i^-_x)^\ast (\omega)$ and $h_x(0)=0.$

There exists $\rho' >0$ so that :
\begin{enumerate}
\item for any $t$ with $t >\rho',$ 
$x\in \mathcal X$  and $a \in \Omega^{\ind x}(M)$ the  integral
$$ \int_{\mathbb R^{\ind x}}  e ^{-th_x}(i^-_x)^\ast (a)\in \mathbb C$$
is absolutely convergent, 
\item   the map $a \rightsquigarrow  \int_{\mathbb R^{\ind x}}  e ^{-th_x}(i^-_x)^\ast (a)\in \mathbb C$ defines  the  linear maps 
$Int^k:  \Omega^k(M) \to C^k(M;X)$  which, when restricted to $(\Omega^k(M)_{\text {sm}} ,$ provide an isomorphism from  
$(\Omega^\ast(M)_{\text {sm}} (M)(t), d^\omega (t))$ to $\mathbb C^\ast (M, X, \omega, \mathcal O)(t).$
\end{enumerate}
\end{theorem}

Theorem \ref {th3} identifies  the cochain complex $\mathbb C^\ast (M, X, \omega, \mathcal O)(t),$  an object determined by dynamics,  to  a subcomplex of $(\Omega^\ast(M) , d^\ast_\omega(t)), $ precisely to $(\Omega^\ast(M)_{\text {sm}} (M)(t), d^\omega (t)).$
Moreover  $(\Omega^\ast(M)_{\text {sm}} (M)(t), d^\omega (t))$ gets a canonical base, the image by $(Int^\ast )^{-1}$ of the base of $C^\ast (M; X)$ 
determined by  the rest points. Note that this  base depends only on the integration  on unstable manifolds,  the metric $g$ and the closed one form  $\omega$ (which ultimately determines $\Omega^\ast _{sm}(M) (t)$).   However, if one regards $d^\ast_\omega(t)$  
with respect to this base as a matrix, its entries are exactly the functions $I^{\mathcal O, \omega}_{x,y} (t)$ which count the  instantons.  

Theorem \ref{th3} is first proved for $X= -\grad_g\omega.$ To obtain the result as stated one needs 
additional arguments.
There is 
a qualitative difference between these two vector fields.
At  rest points the linearization of the first vector field has all eigenvalues real numbers  $\ne 0$  while for  the second vector field the eigenvalues can be complex numbers only with real part  $\ne 0.$ 

\vskip .2in
{\bf The invariant $\mathcal R$:}

Consider $(M,g)$ a  Riemannian manifold of dimension $n$ and $\omega$ a closed one form. Let $\Psi (g)\in \Omega^{n-1}(TM\setminus  M; \pi^\ast \mathcal O_M )$ be the global angular form  associated to the Riemannian manifold $(M,g)$  introduced by  Mathai-Quillen  \cite {MQ86}.
 Here $\pi^\ast \mathcal O_M$ denotes the pull back  of the orientation bundle  $\mathcal O_M $ by tangent bundle map on $TM\setminus M.$

Suppose $\omega\in \Omega^1 (M)$ is a real valued closed one form 
and $X:M\to TM$ a smooth vector field with no rest points.  In \cite {BH03} and  \cite  {BH06} the following quantity
 $$\mathcal R(X,g,\omega):=\int _M \omega\wedge X^\ast (\Psi (g))$$ was introduced  as  a numerical invariant of the triple $(X, g, \omega).$ It was noticed that this invariant can be extended to vector fields with isolated rest points (in particular hyperbolic) using a "geometric regularization" of the possibly divergent integral 
$\int _{M\setminus \mathcal X} \omega\wedge X^\ast (\Psi (g)).$ The invariant has a number of remarkable properties which have been presented in \cite  {BH06}.

\vskip .2in

{\bf  The function $\log Vol (t):$}

Regard  ${Int ^k}(t): \Omega^k(M)_{\text {sm}} (t) \to C^k (M, X)$  as an isomorphism between  two finite dimensional vector spaces  equipped with scalar products\footnote {If the vector space is over $\mathbb C$ then "scalar product" means "Hermitian scalar product".}.
The first with the scalar product induced from the Riemannian metric $g,$ the second  with  the only scalar product which makes the base provided by the rest points orthonormal.  Recall that for $\alpha: V\to W$ an isomorphism between two vector spaces with scalar products 
\begin{equation*}
\log \rm{Vol} (\alpha)= 1/2 \log \det (\alpha^\sharp\cdot \alpha)
\end{equation*} 
with $\alpha^\sharp$ the adjoint of $\alpha.$
Write 
 $\log \rm{Vol }_k(t):= \log\rm{Vol} (Int^k(t))$
 and define
\begin{equation}\label {E2}
\log\rm{Vol}(t):= \sum (-1)^k \log \rm{Vol}_k(t).
\end{equation}

{\bf The Ray- Singer large torsion $\log T_{la}(t):$}

Recall that for a positive  self adjoint elliptic pseudo differential  operator $\Delta,$  the  spectrum  is a countable collection of positive real 
 numbers $0 < \lambda_1 \leq \lambda_2 \leq  \cdots \leq \lambda_k \leq \lambda_{k+1} \cdots  .$   
 Ray and Singer have defined the regularized determinant  
 $\det \Delta$ by  the formula

\begin{equation}\label{E3}
\log \det \Delta = -d/ds_{s=0} \zeta^{\Delta}(s)
\end{equation}
 where $\zeta^\Delta(s),$  the zeta function of $\Delta,$ is the analytic continuation of  $\sum_i ( \lambda_i)^{-s}$ well defined  for $\Re s >0.$ It is known that  $\zeta^\Delta (s)$ is  a  meromorphic function in  the entire complex plane and has $0$ as a regular value.
The definition  can be applied to the  operator $\Delta^{\omega, \text{la}}_k(t).$   Define 
 
 \begin{equation}\label {E4}
 \log T_{\text{ la}}(t):=  1/2 \sum _k k (-1)^{k+1} \log\det \Delta^{\omega, \text{la}}_k(t)
 \end{equation}
 which is a real analytic function\footnote {Actually $e^{2 \log T_{\rm{la}}(t)}$  can be extended to a  holomorphic function in $z$ for $\Re z>\rho$.} in $t$ for $t>\rho$

\begin{theorem}\label{th4}
Suppose $\omega$ is  Morse Lyapunov for $X,$  and  $X$  satisfies {\bf G} and {\bf EG} and $g$ is a Riemannian metric flat near the rest points of $X$.  Then there exists $\rho>0$ so that for $t>\rho$   

\begin{equation}\label {E5}
\log T_{\text{la}} (M, g, \omega) (t) - \log \rm{Vol}(t) +  t\mathcal R( g, X, \omega)
=  Z^X_{[\omega]} (t).
\end{equation}
\end{theorem}
The right side of the above equality is a dynamical quantity while the left side in a spectral geometry quantity.

Next observe that  the family of operators $$\Delta^\omega_q(z)= \Delta_q + z (L_X + L_X^\ast) + z^2||\grad_g \omega||^2$$ is a selfadjoint holomorphic family  of type A in the sense of Kato  \cite{K76}, and therefore there exists a family  of functions $\lambda_n(z)$ and a family of differential  forms $a_n (z)\in \Omega(M)^q$ both holomorphic in $z$ in a  neighborhood of $[0,\infty)$ in the complex plane so that 
 \begin{enumerate}
 \item 
 $\lambda_n(z)'$s  exhaust the eigenvalues of $\Delta^\omega_q(z),$
 \item 
 $a_n(z)$ is an eigenform corresponding to the eigenvalue $\lambda_n(z).$
 \end{enumerate} 

As a consequence the left side of  equality (\ref{E5}) has an analytic continuation to a neighborhood of $[0,\infty)$ and since the right side is a well defined holomorphic function for $\Re z > \rho$  so is the left side.   
It also follows from Kato's theory that:

\begin{theorem}\label{th5}
The two  complexes   $(\Omega^\ast(M)_{\text {sm}} (M)(t), d^\omega (t))$ and $(\Omega^\ast(M)_{\text {la}} (M)(t), d^\omega (t))$
have analytic continuation
to a neighborhood of $[0,\infty).$
\end{theorem}
One expects this to be the case for  $\Re z>0.$ 
\vskip .1in

The proofs of Theorems \ref {th3} and \ref {th4} are  contained in  \cite{BH04}  and can be also derived from \cite {BH03}
while Theorem \ref{th5} is a consequence of Kato's theory. 

Theorems  \ref {th3},  \ref {th4} and \ref{th5} can be generalized. Both the deRham complex $(\Omega^\ast (M), d^\ast)$ and the Witten deformation 
$\Omega^\ast (M), d^\ast_\omega(t)= d^\ast +t\omega \wedge$ can be twisted by a closed complex valued one-form $\eta.$ Precisely $(\Omega^\ast (M), d^\ast)$ can be replaced by 
$(\Omega^\ast (M), d^\ast_\eta= d^\ast + \eta \wedge )$ and $(\Omega^\ast (M), d^\ast_\omega(t)= d^\ast +t\omega \wedge)$
by $(\Omega^\ast (M), d^\ast_{\eta; \omega}= d^\ast + (\eta+t\omega) \wedge ).$

Theorems  \ref {th3},  \ref {th4} and \ref{th5}  remain true with proper modification. For example the conclusion  of  Theorem  \ref{th4} remains the same  by replacing $t\omega$ by $\eta +t\omega$ and $Z^{[\omega]} (z)$ by $Z^{[\eta]; [\omega]} (z)$
with 

$$ Z^{\xi_1. \xi_2}(z)=: \sum _{\gamma \in [S^1, M]} \mathcal Z (\gamma)e^{-(\xi_1 +z\xi_2) (\gamma)} $$
for $\xi_1\in H^1(M; \mathbb C), \xi_{2}\in H^1 (M;\mathbb R).$
An equivalent form of  this stronger  result is proved  in \cite {BH03}.

{ \bf An interesting particular case:}

If $X$ has no rest points then {\bf G} and {\bf EG} reduce to  {\bf NCT} i.e.   all closed trajectories are nondegenerate. 
Moreover $(\Omega^\ast (M), d^\ast_\omega (t))  =(\Omega^\ast_{\rm{la}} (t), d^\ast_\omega(t))$  and for any $t>0$ the operator $\Delta^\omega_k (t)$ is invertible and therefore positive for any $k.$

The following statement is a minor improvement of a result of J. Marksick 
\cite  { BH04}.

\begin{corollary} \label {C1}Suppose $X$ is a smooth vector field with no rest points  and $\xi\in H^1(M;\mathbb R)$ and $\omega$ a Lyapunov form representing $\xi$
and $g$ a Riemannian metric  on $M.$ Suppose that all closed trajectories are nondegenerate.   Denote by $$\log T_{an}(t):= 1/2 \sum_q (-1)^q \log \det \Delta_q^\omega (t).$$
Then 
$$\log T_{an}(t) + t\int  \omega \wedge X^\ast(\Psi_g) = Z^\xi (t).$$
\end{corollary}
Corollary \ref {C1} is a particular version of Theorem \ref{th4}.
\section {Topology}

In this section we review  results relating elements of dynamics with topological invariants. 
\vskip .1in
{\bf Twisted cohomology $H^\ast (M;\xi)$:} 

Recall that a cohomology class $\xi\in H^1(M;\mathbb C)$  permits the definition of twisted cohomology $H^\ast (M;\xi).$ This can be defined in any setting, simplicial, singular, Cech, deRham. All settings lead to isomorphic vector spaces at least  for a smooth manifold $M.$
The deRham version of $H^\ast (M;\xi)$ is defined as  the cohomology of the complex of smooth differential forms equipped  with the differential $d^\ast_\eta: \Omega ^\ast \to \Omega^{\ast +1}$ with $d_\eta = d + \eta \wedge.$ Here  $\eta$ is a closed one form representing $\xi$ (in de Rham cohomology). 
The twisted cohomology  for $t[\omega] \in H^1(M:\mathbb R)\subset H^1(M;\mathbb C)$ is therefore the cohomology of the complex $(\Omega^\ast, d^\ast_{\omega}(t))$ in Witten deformation. 
It is not hard to see that given $\xi$ the dimension of $H^\ast (M; z\xi)$ changes for only finitely many $z$ so we write $\underline{\beta}_q^\xi$ for $\dim H^q(M; t\xi)$  for $t$ large.
The following well known result of Novikov provides strong restrictions  for the numbers  of rest points of a vector field which satisfies {\bf H, MS
} and {\bf L}.

\begin{theorem} (Novikov)

If $X$ satisfies {\bf H, MS} and {\bf L} with $\xi$ a Lyapunov  class for $X$ then:
\vskip .1in

\begin{equation*}
\begin{aligned}
n_k &\geq \underline{\beta}^{\xi}_k\\
\sum_{0\leq k\leq r}  (-1)^k n_k &\geq \sum_{0\leq k\leq r}\underline{ \beta}^{\xi}_k , \ r \ \text{even}\\
\sum_{0\leq k\leq r}  (-1)^k n_k &\geq \sum_{0\leq k\leq r} \underline{\beta}^{\xi}_k,  \ r \ \text{odd}.
\end{aligned}
\end{equation*}
\end{theorem}
\vskip .1in

 

{\bf Torsion: }

The cochain complex 
$(C^\ast (M; X), \delta^\ast_{\mathcal O,\omega},(z))$ is equipped with a base and  has well  defined  {\it torsion}   not  explained here in the full generality.
In case that the $z[\omega]$-twisted cohomology  is trivial  the torsion is a  complex number defined up to a sign and  its  square  depends holomorphically  on $z.$
The domain of this function consists of the complex numbers $z$ with $H^\ast (M; z[\omega])=0$ which  is the complement of a finite set  provided that is nonempty.
 \vskip .1in 
The cohomology class $z\xi,$ $\xi\in H^1(M;\mathbb C)$   can be interpreted as a rank one  complex representation of the fundamental group.
The manifold $M$ equipped with  this  representation and  
with  an {\it Euler structure}   has a  {\it Milnor--Turaev torsion}  not described in this paper. 
In case that the $z\xi$-twisted cohomology is trivial the Milnor--Turaev torsion  is a complex number  defined up to sign  and its square  depends holomorphically 
on $z.$ The domain of this function  consists of the complex numbers $z$ with $H^\ast (M; z \xi)=0.$  
The vector field $X$ and some minor additional data determine an Euler structure.

  It is possible to show that  if $X$ satisfies {\bf G, L} and {\bf EG} with $\xi= [\omega]$ a Lyapunov cohomology class  for $X$ then 
 the torsion of $(C^\ast (M; X), \delta^\ast_{\mathcal O,\omega}(z))$ multiplied by  the function $e^{Z^\xi (z)}$ is essentially  the Milnor--Turaev torsion  of $M, z\xi$ and the Euler class defined by $X.$ 
For more details consult \cite {BH07}. This relates the functions $I^{\mathcal O, \omega}_{x,y}(z),$  $Z^\xi (z)$ and the topology.
The statement above is a reformulation and a minor extension of results of  
 Hutchings--Lee and Pajitnov  \cite {HL99} and \cite {P03}  which will not be explained here in details. However,
we will  explain  below this statement in case $t$ is a real number.
\vskip .1in
Suppose $H^\ast (M; t\omega)=0$ for $t>\rho.$ For such $t,$ $\Delta_\ast ^\omega (t)$ and then $\Delta_{\ast} ^{\omega, \rm{sm}},$ are invertible and therefore have non vanishing determinants.
The determinant of $\Delta^\omega_k$  refers to the $\zeta-$ regularized determinant as explained  in the previous section and defined by formula (\ref{E2}). In particular, we have $\log T_{\rm{an}} (t)$ and $\log T_{\rm{sm}}(t) $  defined by the formula (\ref{E3}) in the previous section.  Clearly 

\begin{equation}\label{E6}
\log T_{\rm {an}}(t) =\log T_{\rm{sm}} (t )  +\log T_{\rm {la}} (t).
\end{equation}
\vskip .1in
The cochain complex  $(C^\ast (M; X), \delta^\ast_{\mathcal O,\omega},(z))$ has a base and therefore we can regard its components as equipped with the  scalar product which makes the base orthonormal and implicitly define the corresponding Laplacians $\Delta_{X}(t).$ 

Note also that by Theorem \ref {th2} the cohomology of the finite dimensional cochain complex $(C^\ast (M; X), \delta^\ast_{\mathcal O,\omega},(z))$ for $t$ large ($t >\sup \{\rho, \rho'\}$)
 is trivial  and then by  formula  (\ref{E4}) defines $\log T_{X}(t)$  for $t$ large.   
As $\log \rm{Vol}(t) = \log T_{\rm{sm}}(t) - \log T^\omega_{X}(t)$  combining  (\ref {E5}) and (\ref {E6})
one concludes  that :

$$\log T_{\rm{an}}(t) - t \mathcal R(M,\omega, g)= \log T_{X}(t) + Z^{[\omega]}(t) $$ which, in the case $X$  satisfies in addition to  {\bf G} and {\bf L} the property {\bf EG},  is  a generalization of Marsick's  result to the case  when $X$ has rest points.  This statement  is equivalent via the work of Bismut--Zhang 
\cite {BZ92} to an analytic version of Hutchings- Lee Pajitnov theorem.


\end{document}